\documentclass[11pt]{amsart}
\usepackage{amssymb}
\usepackage{amsmath}
\usepackage[active]{srcltx}
\usepackage{t1enc}
\usepackage[latin2]{inputenc}
\usepackage{verbatim}
\usepackage{amsmath,amsfonts,amssymb,amsthm}
\usepackage[mathcal]{eucal}
\usepackage{enumerate}
\usepackage[centertags]{amsmath}
\usepackage{graphicx}
\usepackage{color}
\graphicspath{ {./images/} }

\newtheorem{theorem}{Theorem}

\newtheorem{corollary}{Corollary}

\begin{document}
\author{N. Anakidze, N. Areshidze, L.-E. Persson and G. Tephnadze}

\title[Approximation by $T$ means]{Approximation by $T$ means with respect to  Vilenkin system in Lebesgue spaces and Lipschitz classes}
\address {N. Anakidze, The University of Georgia, School of Science and Technology, 77a Merab Kostava St, Tbilisi, 0171, Georgia.}
\email{nino.anakidze@mail.ru}
\address{N. Areshidze, Tbilisi State University, Faculty of Exact and Natural Sciences, Department of Mathematics, Chavchavadze str. 1, tbilisi 0128, georgia}
\email{nika.areshidze15@gmail.com}
\address {L.-E. Persson, UiT The Arctic University of Norway, P.O. Box 385, N-8505, Narvik, Norway and department of Mathematics and Computer Science, Karlstad University, 65188 Karlstad, Sweden.}
\email{larserik6pers@gmail.com}
\address {G. Tephnadze, The University of Georgia, School of Science and Technology, 77a Merab Kostava St, Tbilisi, 0171, Georgia.}
\email{g.tephnadze@ug.edu.ge}

\thanks{The research was supported by Shota Rustaveli National Science Foundation (SRNSF) grant no. FR-24-698.}

\thanks{}

\begin{abstract}
In this paper we present and prove some new results concerning approximation properties of $T$ means with respect to the Vilenkin system in Lebesgue spaces and Lipschitz classes for any $1\leq p<\infty$. As applications, we obtain extension of some known approximation inequalities.
\end{abstract}
\maketitle 
\date{}



\section{Introduction}

Concerning some definitions and notations used in this introduction we refer to Section 2.   

It is well-known (see e.g. \cite{gol}, \cite{PTWbook} and \cite{Zy}) that,
 for any $1\leq p\leq \infty$ and $f\in L_p(G_m),$ there exists an absolute constant $C_p,$ depending only on $p$ such that
 
\begin{equation*}
\left\Vert \sigma_nf\right\Vert_p\leq C_p\left\|f\right\|_{p}.
\end{equation*}
Moreover, (see e.g. \cite{PTWbook})	if $1\leq p\leq \infty $, $M_N\leq n< M_{N+1}$, $f\in L^p(G_m)$ and $n\in \mathbb{N},$ then

\begin{eqnarray}\label{aaa}
\left\Vert \sigma_{n}f-f\right\Vert_{p}
\leq 2R^5\sum_{s=0}^{N} \frac{M_{s}}{M_{N}}\omega_p\left(1/M_s,f\right),
\end{eqnarray}
where $R:=\sup_{k\in\mathbb{N}}m_k$ and $\omega_p(\delta,f)$ is the modulus of continuity of $L^p$ functions, $1\le p <\infty$  functions defined by

\begin{eqnarray*}
	\omega_p(\delta,f)=\sup_{|t|<\delta}\Vert f(x-t)-f(x) \Vert_p, \quad \quad \delta>0.
\end{eqnarray*}


It follows that if $f\in Lip\left( \alpha ,p\right) ,$ i.e.,
\begin{eqnarray*}
	Lip(\alpha,p):=\{f\in L^p: \omega_p(\delta,f)=O(\delta^{\alpha}) \quad\textrm{as}\quad \delta \to 0 \},
\end{eqnarray*}
then	
	\begin{equation*}
	\left\Vert\sigma_nf-f\right\Vert _{p}=\left\{
	\begin{array}{c}O\left(1/M_N\right), \text{ \ \ \ if \ \ \ }\alpha >1,\\
	O\left(N/M_N\right),\text{ \ \ \ if \ \ \ }\alpha=1,\\
	O\left( 1/M_n^{\alpha}\right) ,\text{ \ \ \ \ \ \ \ if }\alpha <1.\end{array}\right.
	\end{equation*}
Moreover, (see e.g.  \cite{PTWbook}) if $1\leq p< \infty ,$ $f\in L^{p}(G_m)$ and

	\begin{equation*}
	\left\Vert \sigma _{M_n}f-f\right\Vert_{p}=o\left( 1/M_n\right), \ \text{as} \ n\rightarrow \infty,
	\end{equation*}
then $f$  is a constant function.

For the maximal operators of Vilenkin-Fejer means $\sigma ^{*}$, defined by
$$\sigma ^{*}f=\sup_{n\in\mathbb{N}}\vert \sigma_n f\vert$$ 
the weak-$(1,1)$ type inequality  
\begin{equation*}
\left\Vert \sigma^{*}f\right\Vert_{weak-L_1} \leq C\left\|
f\right\| _{1},\text{ \quad }\left( f\in L^1(G_m)\right)
\end{equation*}
can be found in  Schipp
\cite{Sc} for Walsh series and in Pál, Simon \cite{PS} and Weisz \cite{We1}  for bounded Vilenkin
series. Boundedness of the maximal operators of Vilenkin-F\'ejer means of the one- and two-dimensional cases can be found in Fridli \cite{Fridli}, G\'at \cite{gat}, Goginava \cite{Goginava1}, Nagy and Tephnadze \cite{NT1,NT2,NT3,NT4}, Simon \cite{Simon1,Simon2} and Weisz \cite{We3}.

Convergence and summability of  N\"orlund means with respect to Vilenkin systems were studied by Areshidze and Tephnadze\ \cite{AT1}, Blahota and Nagy \cite{BN1}, Blahota, Persson and Tephnadze \cite{BPT1} (see also \cite{BNPT1,BPTW1,BT2,BTT1}),  Fridli, Manchanda and  Siddiqi \cite{FMS}, Goginava \cite{Goginava}, Nagy \cite{na,nag,n} (see also  \cite{BN} and \cite{BNT}) and Memic \cite{Memic}.

 M\'oricz and Siddiqi \cite{Mor} investigated the approximation properties of
some special N\"orlund means of Walsh-Fourier series of $L^{p}$ functions in
norm. In particular, they proved that if $f\in L^p(G),$ $1\leq p\leq \infty,$ $n=2^j+k,$ $1\leq k\leq 2^j \ (n\in \mathbb{N}_+)$ and $(q_k,k\in \mathbb{N})$ is a sequence of non-negative numbers, such that
\begin{equation*}\label{MScond}
\frac{n^{\gamma-1}}{Q_n^{\gamma}}\sum_{k=0}^{n-1}q^{\gamma}_k =O(1),\ \ \text{for some} \ \ 1<\gamma\leq 2,
\end{equation*}
then there exists an absolute constant $C_p,$ depending only on $p$ such that

\begin{equation*} \label{MSes}
\Vert t_nf-f\Vert_p\leq \frac{C_p}{Q_n} \sum_{i=0}^{j-1}2^iq_{n-2^i}\omega_p\left(\frac{1}{2^i},f\right)+C_p\omega_p\left( \frac{1}{2^j},f\right),
\end{equation*}
when the sequence $(q_k,k\in \mathbb{N})$ is non-decreasing, while 

$$
\Vert t_nf-f\Vert_p\leq \frac{C_p}{Q_n} \sum_{i=0}^{j-1}\left(Q_{n-2^i+1}-Q_{n-2^{i+1}+1} \right)\omega_p\left(\frac{1}{2^i},f\right) +C_p\omega_p\left( \frac{1}{2^j},f\right),
$$
when the sequence  $(q_k,k\in \mathbb{N})$ is non-increasing.


Tutberidze \cite{Tut2} (see also \cite{PTWbook}) proved that if $T_n$ are $T$ means generated by either a non-increasing sequence $\{q_k,k\in \mathbb{N}\}$ or a  non-decreasing sequence $\{q_k,k\in \mathbb{N}\}$ satisfying the condition

\begin{equation*}\label{cond3}
\frac{q_{0}}{Q_k}=O\left(\frac{1}{k}\right),\text{ \ \ as \ \ }
k\rightarrow \infty,
\end{equation*} 
then there exists an absolute constant $C,$ such that 

\begin{equation*}\label{eq1}
\left\Vert T^{*}f\right\Vert_{weak-L_1} \leq C\left\|
f\right\| _{1},\text{ \quad }\left( f\in L^1(G_m)\right)
\end{equation*}
holds. From these results follows that if $f\in L^p(G_m),$ where $1\leq p< \infty$ and either the sequence  $\{q_k,k\in \mathbb{N}\}$ is non-increasing, or $\{q_k,k\in \mathbb{N}\}$ is a sequence of  non-decreasing numbers, such that the condition

	\begin{equation}\label{cond2}
\frac{q_{n-1}}{Q_n}=O\left(\frac{1}{n}\right),\text{\ \ as \ \ }n\rightarrow \infty,
\end{equation} 
is fulfilled, then
$$ \underset{n\rightarrow \infty }{\lim }\Vert {{T}_{n}}f-f\Vert_p\to 0, \ \ \ \text{as} \ \ \ n\to \infty.$$

For the Walsh system in \cite{MR} M\'oricz and Rhoades proved that if $f\in L^p,$ where $1\leq p< \infty,$  and ${{T}_{n}}$ are regular $T$  means generated by a non-increasing sequence $\{q_k, \ k\in \mathbb{N}\},$  then, for any $2^N\leq n< 2^{N+1},$  we have the following approximation inequality:

\begin{equation}\label{ineq4}
\Vert T_nf-f\Vert_p\leq\frac{C_p}{Q_n}
\sum_{s=0}^{N-1} 2^{s}q_{2^s}\omega_p\left(1/2^s,f\right)
+C_p\omega_p\left(1/2^N,f\right).	
\end{equation}
In the case when the sequence $\{q_k, \ k\in \mathbb{N}\}$ is non-decreasing and satisfying the condition
\begin{equation}\label{cond4}
\frac{q_{k-1}}{Q_k}=O\left(\frac{1}{k}\right),\text{ \ \ as \ \ }
k\rightarrow \infty,
\end{equation} 
then the following inequality holds:

\begin{eqnarray}\label{ineq5}
\Vert T_nf-f\Vert_p \leq C_p\sum_{j=0}^{N-1}2^{j-N}\omega_p\left(1/2^j,f\right)+C_p\omega_p\left(1/2^N,f\right).
\end{eqnarray}

In this paper we use a new approach and generalize the inequalities in  \eqref{ineq4} and \eqref{ineq5} for  $T$ means with respect to the Vilenkin system (see Theorems 1 and 2). We also prove a new inequality for the subsequences $\{T_{M_n}\}$  means when the sequence $\{q_k, \ k\in \mathbb{N}\}$ is non-decreasing (see Theorem 3). 

The main results are presented in Section 3 together with 4 Corollaries, which also generalize known results (see \cite{MR}). The detailed proofs are presented in Section 4. In order not do disturb our presentations in these Sections we have reserved Section 2 for some Preliminaries.

\section{Preliminaries}

Let  $\mathbb{N}_{+}$ denote the set of the positive integers, $\mathbb{N}:=
\mathbb{N}_{+}\cup \{0\}.$ Let 
$m=:(m_0,m_1,...)$
be a sequence of positive integers not less than 2. Denote by
\begin{equation*}
Z_{m_k}:=\{0,1,...,m_k-1\}
\end{equation*}
the additive group of integers modulo $m_k.$ Define the group $G_m$ as the complete direct product of the group $%
Z_{m_k}$ with the product of the discrete topologies of $Z_{m_k}\textrm{'s}$.

The direct product $\mu $ of the measures 
\begin{equation*}
\mu_k\left( \{j\}\right):=1/m_k\text{ \ }(j\in Z_{m_k})  
\end{equation*}
is the Haar measure on $G_m$ with $\mu \left( G_m\right) =1.$

If $\sup_{k\in\mathbb{N}}m_k<+\infty$, then we call $G_m$ a bounded Vilenkin group. If the sequence $\{m_k\}_{k\ge0}$  is unbounded, then $G_m$ is said to be unbounded Vilenkin group. In this paper we consider only bounded Vilenkin groups.

The elements of $G_m$ are represented by the sequences 

\begin{equation*}
x:=(x_{0},x_{1},\dots,x_{k},\dots)\qquad \left( \text{ }x_{k}\in
Z_{m_k}\right).
\end{equation*}
It is easy to give a base for the neighborhood of $G_m$, namely

\begin{eqnarray*}
I_{0}\left( x\right):=G_m, \ \ \
I_{n}(x):=\{y\in G_m\mid y_{0}=x_{0},\dots,y_{n-1}=x_{n-1}\}\text{ }(x\in
G_m,\text{ }n\in \mathbb{N}).
\end{eqnarray*}
For simplicity we also define  $I_n:=I_n(0).$

Next, we define a generalized number system based on $m$ in the following way:
\begin{equation*}
    M_0=:1,\quad M_{k+1}=:m_{k}M_k\quad (k\in\mathbb{N})
\end{equation*}
Then every $n\in \mathbb{N}$ can be uniquely expressed as

$$
n=\sum_{k=0}^{\infty }n_{j}M_j, \ \ \ \text{where } \ \ \ n_{j}\in Z_{m_j}  \ \ \ (j\in \mathbb{
N})$$ 
and only a finite number of $n_{j}`$s differ from zero. Let 

 $$|n|=:\max\{j\in\mathbb{N}, n_j\neq0\}.$$
Moreover, Vilenkin (see \cite{Vi,Vi1,Vi2}) investigated the group $G_m$ and introduced the Vilenkin system $\{{\psi}_j\}_{j=0}^{\infty}$ as
\begin{equation*}
\psi _{n}\left( x\right):=\prod_{k=0}^{\infty }r_{k}^{n_{k}}\left( x\right) 
\text{ \quad }\left( n\in \mathbb{N}\right).
\end{equation*}
where $r_k(x)$ are the generalized Rademacher functions defined by
\begin{equation*}
r_{k}( x):=\exp(2\pi ix_k/m_k), \text{ \ \ }\left(k\in \mathbb{N}\right).
\end{equation*}

These systems include as a special case the Walsh system when $m_k=2$ for any $k\in\mathbb{N}$.

The norms (or quasi-norms) $\left\Vert f\right\Vert _{p},$ $0<p<\infty,$ of the Lebesgue spaces $L^{p}(G_{m})$  are defined by 
\begin{equation*}
\left\Vert f\right\Vert _{p}^{p}:=\int_{G_{m}}\left\vert f\right\vert
^{p}d\mu.
\end{equation*}%

The Vilenkin system is orthonormal and complete in $L^{2}\left( G_m\right)
\,$ (see e.g. \cite{AVD} and \cite{sws}).

If $f\in L^{1}\left(G_m\right) $, we can define the Fourier
coefficients, the partial sums of the Fourier series, the Fejér means, the Dirichlet and Fejér kernels 
with respect to the Vilenkin system in
the usual manner:

\begin{eqnarray*}
	\widehat{f}\left( k\right) &:&=\int_{G_m}f\overline{\psi}_{k}d\mu ,\,%
	\text{\quad }\left(k\in \mathbb{N}\right) , \\
	S_{n}f &:&=\sum_{k=0}^{n-1}\widehat{f}\left( k\right) \psi_{k},\text{
		\quad }\left(n\in \mathbb{N}_{+},\text{ }S_{0}f:=0\right) , \\
	\sigma _{n}f &:&=\frac{1}{n}\sum_{k=1}^{n}S_{k}f,\text{ \quad \ \  }\left(n\in \mathbb{N}_{+}\right).\\
 D_{n}&:&=\sum_{k=0}^{n-1}\psi_{k},\text{ \quad \ \  }\left(n\in \mathbb{N}_{+}\right).\\
 K_n&:&=\frac{1}{n}\sum_{k=1}^{n}D_{k},\text{ \quad \ \  }\left(n\in \mathbb{N}_{+}\right).\
\end{eqnarray*}

Recall that  (see e.g. \cite{AVD} and \cite{PTWbook}),

\begin{equation}\label{dn2.3}
\quad \hspace*{0in}D_{M_{n}}\left( x\right) =\left\{
\begin{array}{l}
\text{ }M_n,\text{\thinspace \thinspace \thinspace  if\thinspace
	\thinspace }x\in I_{n}, \\
\text{ }0,\text{\thinspace \thinspace \thinspace \thinspace \thinspace \thinspace \thinspace \thinspace if
	\thinspace \thinspace }x\notin I_{n},
\end{array}
\right.  
\end{equation}

\begin{eqnarray}\label{dn2.4}
D_{M_n-j}\left( x \right)&=&D_{M_n}\left( x \right)-\overline{\psi}_{M_n-1}(-x)D_{j}(-x)\\
&=&D_{M_n}\left( x \right)-\psi_{M_n-1}(x)\overline{D_{j}}(x),\,\,\ 0\le j<M_n.\nonumber
\end{eqnarray}
\begin{equation} \label{fn5}
n\left\vert K_n\right\vert\leq
2R^2\sum_{l=0}^{\vert n\vert } M_l \left\vert K_{M_l} \right\vert,
\end{equation}
 and

\begin{eqnarray} \label{fn40}
\int_{G_m} K_n (x)d\mu(x)=1,  \ \ \ \ \ 
\sup_{n\in\mathbb{N}}\int_{G_m}\left\vert K_n(x)\right\vert d\mu(x)\leq R^5.
\end{eqnarray}
where $R:=\sup_{k\in\mathbb{N}}m_k.$ Moreover, if $n>t,$ $t,n\in \mathbb{N},$ then 

\begin{equation}\label{lemma2}
K_{M_n}\left(x\right)=\left\{ \begin{array}{ll}
\frac{M_t}{1-r_t(x)},& x\in I_t\backslash I_{t+1},\quad x-x_te_t\in I_n, \\
\frac{M_n+1}{2}, & x\in I_n, \\
0, & \text{otherwise.} \end{array} \right.
\end{equation}

The $n$-th N\"orlund mean $t_n$ and $T$  mean $T_n$  of $f\in L^1(G_m)$ are defined by
\begin{equation*} 
t_nf:=\frac{1}{Q_n}\overset{n}{\underset{k=1}{\sum }}q_{n-k}S_kf
\end{equation*}
	and
\begin{equation*} 
T_nf:=\frac{1}{Q_n}\overset{n-1}{\underset{k=0}{\sum }}q_{k}S_kf, 
\end{equation*}
where 
$$
Q_n:=\sum_{k=0}^{n-1}q_k.$$

Here $\{q_k,\ k\geq 0\}$ is a sequence of nonnegative numbers, where $q_0>0$ and
\begin{equation} \label{1.3}
\lim_{n\rightarrow \infty }Q_{n}=\infty .
\end{equation} 
Then, a $T$ mean generated by $\{q_k,\ k\geq 0\}$ is regular if and only if  the condition (\ref{1.3}) is fulfilled (see \cite{PTWbook}).

It is evident that

\begin{equation*}
T_nf\left(x\right)=\underset{G_m}{\int}f\left(t\right)F_{n}\left(x-y\right) d\mu\left(y\right)
\end{equation*}
where 
\begin{equation}\label{1.3T}
F_{n}:=\frac{1}{Q_n}\overset{n-1}{\underset{k=0}{\sum }}q_{k}D_k,
\end{equation}
which are called the kernels of the $T$  means.

By applying Abel transformation we get the following two useful identities: 

\begin{eqnarray}  \label{2b}
Q_n:=\overset{n-1}{\underset{k=0}{\sum}}q_k\cdot 1=\overset{n-2}{\underset{k=0}{\sum}}(q_k-q_{k+1})k
+q_{n-1}(n-1)
\end{eqnarray}
and

\begin{equation}  \label{2bbb}
T_nf=\frac{1}{Q_n}\left(\overset{n-2}{\underset{k=0}{\sum }}(q_{k}-q_{k+1})k\sigma_{k}f +q_{n-1}(n-1)\sigma_{n-1}f\right).
\end{equation}

\section{The Main Results}

Our first main result reads:
\begin{theorem}\label{Corollary3nnconv} 
	Let $f\in L^p(G_m),$ where $1\leq p< \infty$  and ${{T}_{n}}$ are  $T$  means generated by a non-increasing sequence $\{q_k, \ k\in \mathbb{N}\}.$ 
	 Then, for any $n,N\in \mathbb{N},$ $M_N\leq n< M_{N+1},$  we have the following inequality:
		
	\begin{eqnarray}\label{A} \ \ \
	\Vert T_nf-f\Vert_p\leq\frac{6R^6}{Q_n}
	\sum_{j=0}^{N-1} M_{j}q_{_{M_j}}\omega_p\left(1/M_j,f\right)
	+4R^6\omega_p\left(1/M_N,f\right).	
	\end{eqnarray}
\end{theorem}

Next we state and prove a similar inequality for non-decreasing sequences but under some restrictions.

\begin{theorem}\label{Corollary3nnconv1} 
	Let $f\in L^p(G_m),$ where $1\leq p< \infty$ and  $T_n$ are regular $T$  means generated by a non-decreasing sequence $\{q_k, \ k\in \mathbb{N}\}$. 
	Then, for any $n,N\in\mathbb{N},$ $M_N\leq n< M_{N+1},$ we have the following inequality:
	
	\begin{equation}\label{22} 
	\Vert T_nf-f\Vert_p \leq \frac{6R^6q_{n-1}}{Q_n}\sum_{j=0}^{N-1}M_j\omega_p\left(1/M_j,f\right)+\frac{4R^6q_{n-1}M_N}{Q_n}\omega_p\left(1/M_N,f\right).
	\end{equation}
	If, in addition, the sequence $\{q_k, \ k\in \mathbb{N}\}$ satisfies the condition
    \eqref{cond2}, then the inequality	
	
	\begin{eqnarray}\label{Cond0}
	\Vert T_nf-f\Vert_p \leq C_p\sum_{j=0}^{N}\frac{M_j}{M_N}\omega_p\left(1/2^j,f\right)
	\end{eqnarray}
	holds for some constant $C_p$ only depending on $p.$
\end{theorem}

Finally, we state and prove the third main result for the non-decreasing sequences, where we prove a more precise result then that in \eqref{Cond0} and without the restriction \eqref{cond2} but only for subsequences.
\begin{theorem}\label{Corollary3nnconv0} 
	Let $f\in L^p(G_m),$ where $1\leq p< \infty $ and  ${{T}_{k}}$ are regular $T$ means generated by a non-decreasing sequence $\{q_k, \ k\in \mathbb{N}\}$.
	Then, for any $n\in \mathbb{N},$ the following inequality holds:
	
	\begin{eqnarray}\label{inequal} 
&&\Vert T_{M_n}f-f\Vert_p\\ \notag
&\leq&R^2\sum_{j=0}^{n-1}\frac{M_j}{M_n}\omega _p\left(1/M_j,f\right) \\ \notag
&+&\frac{2R^4}{q_{0}}\overset{n-1}{\underset{j=0}{\sum}}\frac{(n-j)(q_{_{M_n-M_j}})M_j}{M_{n}}\omega _p\left(1/M_j,f\right)\\ \notag
&+&\left(R^2+2\right)\omega _p\left(1/M_n,f\right).
	\end{eqnarray}
\end{theorem}

We also point out the following generalizations of some results in \cite{MR} (in this paper only the Walsh system was considered):
\begin{corollary}
	Let $\{ q_k, \  k \geq 0\}$ be a sequence of non-negative and non-increasing numbers, while in case when the sequence is non-decreasing it is assumed that also the condition \eqref{cond2} is satisfied. If $f\in Lip(\alpha, p)$ for some $\alpha > 0$ and $1 \leq p< \infty,$ then
	
	\begin{equation*}
	\Vert T_nf-f\Vert_p=\left\{
	\begin{array}{l}
	\text{ }O(n^{-\alpha}),\text{\qquad \qquad  if\qquad  }0<\alpha<1, \\
	\text{ }O(n^{-1}\log n),\text{\qquad  if
		\quad \ }\alpha=1,\\
	\text{ }O(n^{-1}),\text{\qquad \qquad   if \quad \ \
	}\alpha>1,
	\end{array}
	\right.  
	\end{equation*}
\end{corollary}

\begin{corollary}
	Let $\{ q_k, \  k \geq 0\}$ be a sequence of non-negative and non-increasing numbers such that 
	
	\begin{eqnarray*}\label{dn2.300}
		q_k\asymp k^{-\beta} \ \ \ \text{for some} \ \ \ 0<\beta\leq 1
	\end{eqnarray*}
	is satisfied. 
	
	If $f\in Lip(\alpha, p)$ for some $\alpha > 0$ and $1 \leq p<\infty,$ then
	
	\begin{eqnarray*}
		\Vert T_nf-f\Vert_p
		=\left\{
		\begin{array}{l}
			\text{ }O(n^{-\alpha}),\text{\qquad \qquad \qquad  \ \ \ \ \    if \ \ \  \ }\alpha+\beta<1, \\
			\text{ }O(n^{-(1-\beta)}\log n+n^{-\alpha}),\text{\ \ \   if 
				\quad \ } \alpha+\beta=1,\\
			\text{ }O(n^{-(1-\beta)}),\text{\qquad \qquad \qquad     if \quad \ \ 
			} \alpha+\beta>1, \ \beta>1,\\
			\text{ }O((\log n)^{-1}),\text{\qquad \qquad \ \ \ \     if \quad \ 
			}\ \beta=1.\\
		\end{array}
		\right.  
	\end{eqnarray*}
\end{corollary}
\begin{corollary}
	Let $\{ q_k, \  k \geq 0\}$ be a sequence of non-negative and non-increasing numbers such
	that 
	
	\begin{eqnarray*}
		q_k\asymp {(\log k)}^{-\beta} \ \ \ \text{for some} \ \ \ \beta>0
	\end{eqnarray*}
	is satisfied. 
	
	If $f\in Lip(\alpha, p)$ for some $\alpha > 0$ and $1 \leq p<\infty,$ then
	
	\begin{equation*}
	\Vert T_nf-f\Vert_p=\left\{
	\begin{array}{l}
	\text{ }O(n^{-\alpha}),\text{\qquad \qquad \qquad \ \ \  if \qquad  }0<\alpha<1, \ \beta>0,\\
	\text{ }O(n^{-1}\log n),\text{\qquad \qquad \ \ \  if  \qquad  
	}\alpha=1,\ 0<\beta<1,\\
	\text{ }O(n^{-1}\log n \log\log n),\text{\ \ \ \  if \quad \ \
	}\ \alpha=\beta=1,\\
	\text{ }O(n^{-1}(\log n)^{\beta}),\text{\qquad \qquad   if \quad 
	}\quad \alpha>1, \ \beta>0. \\
	\end{array}
	\right.  
	\end{equation*}
\end{corollary}

\begin{corollary}  
	Let $f\in L^p(G_m),$ where $1\leq p< \infty $  and
	$\{ q_k, \  k \geq 0\}$ is a sequence of non-negative and non-increasing numbers, while in case when the sequence is non-decreasing it is also assumed that the condition \eqref{cond2} is satisfied. Then, 
	
	$$ \underset{n\rightarrow \infty }{\lim }\Vert {{T}_{n}}f-f\Vert_p\to 0, \ \ \ \text{as} \ \ \ n\to \infty.$$
\end{corollary}

\section{Proofs}

\begin{proof}[Proof of Theorem 1.]
Let $M_N\leq n<M_{N+1}.$ Since  $T_n$ are regular $T$  means generated by the sequence of non-increasing numbers $\{q_k:k\in \mathbb{N}\},$  we can combine \eqref{2b} and \eqref{2bbb} and conclude that
	\begin{eqnarray}\label{B}
&&\Vert T_nf-f\Vert_p \\ \notag
&\leq&\frac{1}{Q_n}\left(\overset{n-2}{\underset{j=0}{\sum }}(q_{j}-q_{j+1})j\Vert \sigma_jf-f\Vert_p +q_{n-1}(n-1)\Vert \sigma_{n-1}f-f\Vert_p\right)\\ \notag
&:=&I+II.
\end{eqnarray}	
Moreover,
\begin{eqnarray}\label{I}
I&=&\frac{1}{Q_n}\overset{M_N-1}{\underset{j=1}{\sum}}\left(q_{j}-q_{j+1}\right)j\Vert\sigma_jf-f\Vert_p\\ \notag
&+&\frac{1}{Q_n}\overset{n-1}{\underset{j=M_N}{\sum}}\left(q_{j}-q_{j+1}\right)j\Vert\sigma_jf-f\Vert_p\\ \notag
&:=&I_1+I_2.
\end{eqnarray}
Now we estimate all terms separately. By applying estimate \eqref{aaa} for $I_1$ we obtain that

\begin{eqnarray}\label{I_1}
I_1&\leq&\frac{2R^5}{Q_n}\overset{N-1}{\underset{k=0}{\sum}}\overset{M_{k+1}-1}{\underset{j=M_k}{\sum}}\left(q_{j}-q_{j+1}\right)j  
\sum_{s=0}^{k} \frac{M_{s}}{M_{k}}\omega_p\left(1/M_s,f\right)\\ \nonumber
&\leq&\frac{2R^6}{Q_n}\overset{N-1}{\underset{k=0}{\sum}}M_{k}\overset{M_{k+1}-1}{\underset{j=M_k}{\sum}}\left(q_{j}-q_{j+1}\right) 
\sum_{s=0}^{k} \frac{M_{s}}{M_{k}}\omega_p\left(1/M_s,f\right)\\\nonumber
&\leq&\frac{2R^6}{Q_n}\overset{N-1}{\underset{k=0}{\sum}}\left(q_{_{M_k}}-q_{_{M_{k+1}}}\right) \sum_{s=0}^{k} M_{s}\omega_p\left(1/M_s,f\right)\\\nonumber
&\leq&\frac{2R^6}{Q_n}\overset{N-1}{\underset{s=0}{\sum}}M_{s}\omega_p\left(1/M_s,f\right)\sum_{k=s}^{N-1} \left(q_{_{{M_k}}}-q_{_{M_{k+1}}}\right)\\\nonumber
&\leq&\frac{2R^6}{Q_n}\overset{N-1}{\underset{s=0}{\sum}}M_{s}q_{_{M_s}}\omega_p\left(1/M_s,f\right).
\end{eqnarray}

Moreover,

\begin{eqnarray}\label{I_2}
	\quad \quad I_2&\leq&\frac{2R^5}{Q_n}\overset{n-1}{\underset{j=M_N}{\sum}}\left(q_{j}-q_{j+1}\right)j  
	\sum_{s=0}^{N} \frac{M_{s}}{M_{N}}\omega_p\left(1/M_s,f\right)\\\nonumber
&\leq&\frac{2R^6M_{N}}{Q_n}\overset{n-1}{\underset{j=M_N}{\sum}}\left(q_{j}-q_{j+1}\right)
\sum_{s=0}^{N} \frac{M_{s}}{M_{N}}\omega_p\left(1/M_s,f\right)\\\nonumber
&\leq&\frac{2R^6q_{_{M_N}}}{Q_n}
\sum_{s=0}^{N}M_{s}\omega_p\left(1/M_s,f\right)\\ \notag
&\leq&\frac{2R^6}{Q_n}\overset{N}{\underset{s=0}{\sum}}M_{s}q_{_{M_s}}\omega_p\left(1/M_s,f\right)\\ \notag
&\leq&\frac{2R^6}{Q_n}\overset{N-1}{\underset{s=0}{\sum}}M_{s}q_{_{M_s}}\omega_p\left(1/M_s,f\right)+2R^6\omega_p\left(1/M_s,f\right).
\end{eqnarray}

For $II$ we have that
\begin{eqnarray}\label{II}
	II&\leq&\frac{2R^5M_{N+1}q_{n-1}}{Q_n}  
	\sum_{s=0}^{N} \frac{M_{s}}{M_{N}}\omega_p\left(1/M_s,f\right)\\ \notag
	&\leq& \frac{2R^6}{Q_n}  
	\sum_{s=0}^{N-1}M_{s}q_{_{M_s}}\omega_p\left(1/M_s,f\right)+2R^6\omega_p\left(1/M_N,f\right).
\end{eqnarray}
	The proof of \eqref{A} is complete by just combining \eqref{B}-\eqref{II}.
	
\end{proof}

\begin{proof} [Proof of Theorem 2.]
	Let $M_N\leq n<M_{N+1}.$ Since  $T_n$ are regular $T$  means, generated by a sequence of non-decreasing numbers $\{q_k:k\in \mathbb{N}\},$  by combining \eqref{2b} and \eqref{2bbb}, we find that		
	\begin{eqnarray} \label{T2}
		&&\Vert T_nf-f\Vert_p \\ \notag	
		&\leq&\frac{1}{Q_n}\left(\overset{n-1}{\underset{j=1}{\sum}}\left(q_{j+1}-q_{j}\right)j\Vert\sigma_jf-f\Vert_p+q_{n-1}(n-1)\Vert\sigma_nf-f
		\Vert_p\right)\\ \notag	
		&:=&I+II.
	\end{eqnarray}	
	Furthermore,	
	\begin{eqnarray}\label{I2}
		I&=&\frac{1}{Q_n}\overset{M_{N}-1}{\underset{j=1}{\sum}}\left(q_{j+1}-q_{j}\right)j\Vert\sigma_jf-f\Vert_p\\ \notag
		&+& \frac{1}{Q_n}\overset{n-1}{\underset{j=M_{N}}{\sum}}\left(q_{j+1}-q_{j}\right)j\Vert\sigma_jf-f\Vert_p\\ \notag
		&=&I_1+I_2.
	\end{eqnarray}
	Analogously to \eqref{I_1} we get that
	
	\begin{eqnarray}\label{I12}
		I_1&\leq&\frac{2R^6}{Q_n}\overset{N-1}{\underset{k=0}{\sum}}\left(q_{M_{k+1}}-q_{M_{k}}\right)
		\sum_{s=0}^{k} M_{s}\omega_p\left(1/M_s,f\right)\\ \notag
		&\leq&\frac{2R^6}{Q_n}\sum_{s=0}^{N-1}M_{s}	\omega_p\left(1/M_s,f\right)\overset{N-1}{\underset{k=s}{\sum}}\left(q_{_{M_{k+1}}}-q_{_{M_{k}}}\right)\\ \notag
		&=&\frac{2R^6}{Q_n}\sum_{s=0}^{N-1}M_{s}	\omega_p\left(1/M_s,f\right)(q_{_{M_{N}}}-q_{_{M_s}})\\ \notag
		&\leq&\frac{2R^6q_{M_{N}}}{Q_n}\sum_{s=0}^{N-1}M_{s}	\omega_p\left(1/M_s,f\right)\\ \notag
		&\leq&\frac{2R^6q_{n-1}}{Q_n}\sum_{s=0}^{N-1}M_{s}	\omega_p\left(1/M_s,f\right).
	\end{eqnarray}

	
	In the similar way as in \eqref{I_2} we find that
	\begin{eqnarray}\label{I22}
		I_2&\leq&\frac{2R^5}{Q_n}\overset{n-1}{\underset{j=1}{\sum}}\left(q_{j+1}-q_{j}\right)j  \sum_{s=0}^{N} \frac{M_{s}}{M_{N}}\omega_p\left(1/M_s,f\right)\\ \notag
		&=&\frac{2R^5}{Q_n}\left((n-1)q_{n-1}-Q_{n}\right)\sum_{s=0}^{N} \frac{M_{s}}{M_{N}}\omega_p\left(1/M_s,f\right)\\ \notag
		&\leq &\frac{2R^5M_{N+1}q_{n-1}}{Q_nM_N}\sum_{s=0}^{N} M_{s}\omega_p\left(1/M_s,f\right)\\ \notag
		&\leq&\frac{2R^6q_{n-1}}{Q_n}\sum_{s=0}^{N} M_{s}\omega_p\left(1/M_s,f\right)\\ \notag
		&\leq&\frac{2R^6q_{n-1}}{Q_n}\sum_{s=0}^{N} M_{s}\omega_p\left(1/M_s,f\right)+\frac{2R^6q_{n-1}M_N}{Q_n}\left(1/M_N,f\right).
	\end{eqnarray}
	
	For $II$ we have that
	\begin{eqnarray}\label{II2}
		II&\leq&\frac{2R^5q_{n-1}M_{N+1}}{Q_n}\sum_{s=0}^{N} \frac{M_{s}}{M_{N}}\omega_p\left(1/M_s,f\right)\\ \notag
		&\leq& \frac{2R^6q_{n-1}}{Q_n} \sum_{s=0}^{N}M_{s}\omega_p\left(1/M_s,f\right)\\ \notag
		&=&\frac{2R^6q_{n-1}}{Q_n}\sum_{s=0}^{N-1} M_{s}\omega_p\left(1/M_s,f\right)+\frac{2R^6q_{n-1}M_N}{Q_n}\left(1/M_N,f\right).
	\end{eqnarray}

By combining \eqref{T2}-\eqref{II2} we find that \eqref{22} holds.
Moreover, by using condition \eqref{cond2} we obtain the estimate \eqref{Cond0} so the proof is complete.
\end{proof}

\begin{proof}[Proof of Theorem 3.]
	According to \eqref{dn2.4} we find that
	\begin{eqnarray*} \label{1.21}
	T_{M_n}f=D_{M_n}\ast f-\frac{1}{Q_{M_n}}\overset{M_n-1}{\underset{k=0}{\sum }}q_{k}\left(\left( \psi_{M_n-1}\overline{D_{k}}\right)\ast f\right).
	\end{eqnarray*}
	Hence, by using Abel transformation we get that
	
	\begin{eqnarray*} 	\label{2cc} \ \ \ \ \
	T_{M_n}f&=&D_{M_n}\ast f\\ \notag
	&-&\frac{1}{Q_{M_n}}\overset{M_n-2}{\underset{j=0}{\sum}}\left(q_{_{M_n-j}}-q_{_{M_n-j-1}}\right) j((\psi_{_{M_n-1}}\overline{K_j})\ast f)\\ \notag
	&-& \frac{1}{Q_{M_n}}q_{_{M_n-1}}(M_n-1)(\psi_{_{M_n-1}}\overline{K}_{M_n-1}\ast f)\\ \notag
	&=&D_{M_n}\ast f \\ \notag
	&-&\frac{1}{Q_{M_n}}\overset{M_n-2}{\underset{j=0}{\sum}}\left(q_{_{M_n-j}}-q_{_{M_n-j-1}}\right) j((\psi_{_{M_n-1}}\overline{K_j})\ast f)\\ \notag
	&-& \frac{1}{Q_{M_n}}q_{_{M_n-1}}M_n(\psi_{_{M_n-1}}\overline{K}_{M_n}\ast f)\\ \notag
	&+& \frac{q_{_{M_n-1}}}{Q_{M_n}}(\psi_{_{M_n-1}}\overline{D}_{M_n}\ast f)
	\end{eqnarray*}
	so that
	
	\begin{eqnarray}\label{2c}  
	&&T_{M_n}f(x)-f(x)\\ \notag
	&=&\int_{G_m}(f(x-t)-f(x))D_{M_n}(t)dt\\ \notag
	&-& \frac{1}{Q_{M_n}}\overset{M_n-2}{\underset{j=0}{\sum}}\left(q_{_{M_n-j}}-q_{_{M_n-j-1}}\right) j\int_{G_m}\left(f(x-t)-f(x)\right)\psi_{_{M_n-1}}(t)\overline{K}_j(t)dt\\  \notag
	&-& \frac{1}{Q_{M_n}}q_{_{M_{n}-1}}M_n\int_{G_m}(f(x-t)-f(x))\psi_{_{M_n-1}}(t)\overline{K}_{M_n}(t)dt\\ \notag
	&+& \frac{q_{_{M_{n}-1}}}{Q_{M_n}}\int_{G_m}(f(x-t)-f(x))\psi_{_{M_n-1}}(t)\overline{D}_{M_n}(t)dt\\ \notag
	&=:&I+II+III+IV.
	\end{eqnarray}
	By combining generalized Minkowski's inequality and   \eqref{dn2.3}  we find that
	\begin{eqnarray}\label{fejaprox3}
	\Vert I\Vert _p\leq \int_{I_n}\Vert f(x-t)-f(x))\Vert_p D_{M_n}(t)dt\leq \omega_p\left(1/M_n,f\right).
	\end{eqnarray}
and
\begin{eqnarray}\label{fejaprox4}
	\Vert IV\Vert _p\leq \int_{I_n}\Vert f(x-t)-f(x))\Vert_p {D}_{M_n}(t)dt\leq \omega_p\left(1/M_n,f\right) .
	\end{eqnarray}
Moreover, since
$$ M_nq_{M_n-1}\leq Q_{M_n} , \ \ \text{ for any } \ \ \ n\in \mathbb{N},$$ 
we can use  \eqref{lemma2} and generalized Minkowski's inequality to find that	
	\begin{eqnarray}\label{fejaprox2}
		\Vert III\Vert_p
		&\leq&\int_{G_{m}}\left\Vert f\left( x-t\right) -f\left( x\right) \right\Vert_p
		\left|\overline{K}_{M_n}\left(t\right)\right| d\mu (t)\\ \notag
		&=&\int_{I_{n}}\left\Vert f\left( x-t\right) -f\left( x\right) \right\Vert_p \left|\overline{K}_{M_n}\left(t\right)\right|d\mu(t)\\ \notag
		&+&\sum_{s=0}^{n-1}\sum_{n_s=1}^{m_s-1}\int_{I_{n}\left(n_s e_{s}\right)}\left\Vert f\left( x-t\right) -f\left( x\right) \right\Vert_p \left|\overline{K}_{M_n}\left(t\right)\right| d\mu (t)\\ \notag
		&\leq&\int_{I_{n}}\left\Vert f\left( x-t\right)-f\left( x\right) \right\Vert_p\frac{M_{n}+1}{2}d\mu(t)\\ \notag
		&+&\sum_{s=0}^{n-1}M_{s+1}\sum_{n_s=1}^{m_s-1}\int_{I_n\left(n_se_s\right) }\left\Vert f\left(x-t\right)-f\left(x\right) \right\Vert_pd\mu(t)\\ \notag
		\\ \notag
		&\leq& \omega _{p}\left( 1/M_{n},f\right) \int_{I_{n}}\frac{M_{n}+1}{2}d\mu(t)\\\notag
		&+&\sum_{s=0}^{n-1}M_{s+1}\sum_{n_s=1}^{m_s-1}\int_{I_n\left(n_se_s\right)}\omega _{p}\left( 1/M_s,f\right)d\mu(t)\\ \notag
		&\leq&\omega _{p}\left( 1/M_{n},f\right)+R^2\sum_{s=0}^{n-1}\frac{M_s}{M_n}\omega _p\left(1/M_s,f\right)\\ \notag
		&\leq& R^2\sum_{s=0}^{n}\frac{M_s}{M_n}\omega _p\left(1/M_s,f\right).
	\end{eqnarray}
 	From this inequality and the estimates in \eqref{fejaprox2} it follows also that	
 	
\begin{eqnarray*}\label{fejaprox2sssaaa}
&&M_n\int_{G_m}\left\Vert f\left( x-t\right) -f\left( x\right) \right\Vert_p |\overline{K}_{M_{n}}(t)|\, d\mu (t)\\ \notag
&\leq& R^2\sum_{s=0}^{n}M_s\omega _p\left(1/M_s,f\right).
\end{eqnarray*}

Let $M_k\leq j< M_{k+1}$ By applying  \eqref{fn5} and the last estimate we find that
	
\begin{eqnarray*} \label{fejaprox22a}
        && j\int_{G_m}\left\Vert f\left( x-t\right) -f\left(x\right)\right\Vert_p
        |\overline{K}_{j}(t)|d\mu (t)\\ \notag
        &\leq& 2R^4\sum_{l=0}^{k}\sum_{s=0}^{l}M_s\omega _p\left(1/M_s,f\right).
\end{eqnarray*}

	Hence, by also using \eqref{fn5}  we obtain that	
	\begin{eqnarray}\label{II3}
		\Vert II\Vert_p\\ \notag
		\end{eqnarray}
		\begin{eqnarray*}
		 &\leq&\frac{1}{Q_{M_n}}\overset{M_n-1}{\underset{j=0}{\sum}}\left(q_{_{M_n-j}}-q_{_{M_n-j-1}}\right) j\int_{G_m}\Vert f(x-t)-f(x)\Vert_p \vert\overline{K}_j(t)\vert d\mu (t)\\  \notag
		&\leq&\frac{1}{Q_{M_n}}\overset{n-1}{\underset{k=0}{\sum}}\overset{M_{k+1}-1}{\underset{j=M_k}{\sum}}\left(q_{_{M_n-j}}-q_{_{M_n-j-1}}\right) j\int_{G_m}\Vert f(x-t)-f(x)\Vert_p \vert \overline{K}_j(t)\vert d\mu (t)\\  \notag
		&\leq & \frac{2R^4}{Q_{M_n}}\overset{n-1}{\underset{k=0}{\sum}}\overset{M_{k+1}-1}{\underset{j=M_k}{\sum}}\left(q_{_{M_n-j}}-q_{_{M_n-j-1}}\right)\sum_{l=0}^{k}\sum_{s=0}^{l}M_s\omega _p\left(1/M_s,f\right)
		\\  \notag
		&\leq & \frac{2R^4}{Q_{M_n}}\overset{n-1}{\underset{k=0}{\sum}}\left(q_{_{M_n-M_k}}-q_{_{M_{n}-M_{k+1}}}\right) \sum_{l=0}^{k}\sum_{s=0}^{l}M_s\omega _p\left(1/M_s,f\right)
		\\  \notag
		&\leq & \frac{2R^4}{Q_{M_n}}\overset{n-1}{\underset{l=0}{\sum}}\sum_{k=l}^{n-1}\left(q_{_{M_{n}-M_{k}}}-q_{_{M_{n}-M_{k+1}}}\right)\sum_{s=0}^{l}M_s\omega _p\left(1/M_s,f\right)
		\\  \notag
		&\leq & \frac{2R^4}{Q_{M_n}}\overset{n-1}{\underset{l=0}{\sum}}q_{_{M_n-M_l}}\sum_{s=0}^{l}M_s\omega _p\left(1/M_s,f\right) \\ \notag
		&\leq&  \frac{2R^4}{Q_{M_n}}\overset{n-1}{\underset{s=0}{\sum}}M_s\omega _p\left(1/M_s,f\right)\sum_{l=s}^{n-1}q_{M_n-M_l}
		\\ \notag
		&\leq & \frac{2R^4}{Q_{M_n}}\overset{n-1}{\underset{s=0}{\sum}}M_s\omega _p\left(1/M_s,f\right)q_{_{M_{n}-M_s}}(n-s)\\ \notag
		&\leq& 2R^4\overset{n-1}{\underset{s=0}{\sum}}\frac{(n-s)M_s}{M_n}\frac{q_{_{M_n-M_{s}}}}{q_{0}}\omega _p\left(1/M_s,f\right).
	\end{eqnarray*}

	Finally, by combining 	\eqref{fejaprox3}-\eqref{II3} we obtain \eqref{inequal} so the proof is complete.
\end{proof}

\end{document}